\documentclass{notices}
\usepackage{amsfonts,amssymb,amsmath,amscd,graphicx,color}

\newcommand{\R}{\mathbb{R}}

\title{Fractional derivatives: Fourier, elephants, memory effects, viscoelastic materials and anomalous diffusions}

\date{}

\author{Pablo Ra\'ul Stinga
\affil{Pablo Ra\'ul Stinga is Associate Professor of Mathematics at Iowa State University.
His email address is \texttt{stinga@iastate.edu}.}}

\begin{document}

\maketitle

In the late 17th century, Isaac Newton and Gottfried Wilhelm Leibinz created calculus.
The concept of instantaneous rate of change of a quantity with respect
to another one, that is, the derivative of a function, made a profound impact in science and the society at large.
Differential and integral equations are among the most efficient tools in the mathematical description of
many (but not \emph{all}) natural phenomena. Needless to say, calculus is at the core of a lot of the most advanced 
technological accomplishments of humankind.

Let us consider a function $u$ of a single variable that we will denote by $t$
to represent time. Leibniz introduced the notation
$$\frac{d^nu}{dt^n}$$
for derivatives of order $n\geq0$ of $u$ with respect to $t$, with the understanding
that the derivative of order $n=0$ is just the identity operator, $\frac{d^0u}{dt^0}=u$.

The beginning of \emph{fractional calculus} can be traced back to the origins of classical
calculus itself. On September 30, 1695, Leibniz sent a letter from Hanover to
Marquis de l'Hospital in Paris \cite{MR0141575}. Since Leibniz ``had some extra space left to write'' (see \cite[p.~301]{MR0141575}),
he shared with l'Hospital some remarks from the analogy of integer powers of derivatives $d^n$
and integrals $d^{-n}=\int^n$. 
Leibniz introduced the symbol $d^{\frac{1}{2}}$ to denote
a derivative of \emph{fractional order} $n=1/2$.
He did not give a formal definition but believed that one could express a derivative of fractional order
with an infinite series. He wrote
\textit{``There is the appearance that one day we will come to some
very useful consequences of these paradoxes,''} leaving his seminal thoughts
for l'Hospital to think about. It seems that l'Hospital did not pursue the matter any further.
Leibniz also shared his thoughts with Wallis in 1697, using again the notation $d^{1/2}y$.

We may say that fractional calculus was conceived as a question of \emph{extension of meaning}:
given an object defined for integers, what is its extension to non-integers? 
In fact, we do not need to restrict to \emph{fractions} only. We could also ask about the meaning
of the derivative of orders $\sqrt{2}$ or $\pi$ of a function. In this sense, we use the words
\emph{fractional derivative} to refer to derivatives of non-integer order.

A few other initial ideas were introduced by
Euler (1738) and Lacroix (1820) for fractional derivatives of power functions $t^\beta$
and by Laplace (1812) for functions representable by a specific integral.
In the 1820s, Joseph Fourier was the first one to give a definition that worked for \emph{any} function
(in Fourier's view, \emph{``there is no function $f(x)$, or part of a function, which cannot be
expressed by a trigonometric series,''} see \cite[Par.~418, p.~555]{MR1414430}).
The first known application of fractional calculus appeared shortly afterwards when
Abel found an integral of fractional order in his solution to the tautochrone problem.
Further contributions by many others flourished during the 19th and 20th centuries.
The monumental work \cite{MR1347689} contains a quite detailed historical overview
while a brief history of fractional calculus can be read in \cite{1975}, see also \cite{MR444394}.

Real-world phenomena exhibiting long range-interactions,
memory effects, anomalous diffusions and avalanche-like behaviors are 
very well-described by fractional derivative models. Furthermore,
due to their mathematical features, fractional derivatives are becoming central in image processing.
Nowadays, fractional calculus is one of the most vibrant areas in mathematics that continues
to expand. In fact, much is still to be understood about fractional derivatives from the theoretical and computational points of view.

In this article, we will go back to the widely 
overlooked historical definition of fractional derivative given by Fourier.
We will then show, with modern ideas, how his definition is in fact 
the so-called Marchaud--Weyl fractional derivative. After describing some of
the latest analytical advances in the theory of fractional calculus,
we will present three natural applications to processes with memory effects.
The first one is population growth. For the second one, we will go back to the ideas of Boltzmann, combine them with
engineering experiments, and come up with a model for viscoelastic materials driven by fractional derivatives.
The last application is a fractional model based on the anomalous diffusion behavior observed
in many natural systems.

The literature on theory and applications of fractional calculus is massive.
A Google Scholar search of articles containing the words ``fractional derivative''
gives 120,000 results, with 52,500 of them published in the past 20 years.
It is not the scope of this article to be exhaustive by any means.
Therefore, we will only refer to those works that are directly related with the presentation,
leaving many interesting references out.

\section*{Fourier's fractional derivative}

In 1822, Joseph Fourier was finally allowed to publish his 1807 original research in the form of a
comprehensive monograph entitled \emph{``Th\'eorie Analytique de la Chaleur''} \cite{MR1414430}.
In this work, he introduced the heat equation $\partial_tv=\Delta v$.
To solve it, he created the technique that is nowadays taught to every undergraduate student in mathematics,
physics, engineering and computer science: the method of separation of variables. Fourier studied
the development of an arbitrary function in trigonometric (cosines and sines) series and integrals,
tools that are now known as Fourier series and Fourier transform, respectively.
Towards the end
of his monograph, Fourier looks at representations of real single-variable functions $u$ as
$$u(t)=\frac{1}{2\pi}\int_{-\infty}^\infty dy\,u(y)\int_{-\infty}^\infty d\omega\cos(\omega t-\omega y).$$
(Here we have modified Fourier's notation for the variables to make them consistent with the rest of this article.)
The expression above presents diverse analytical applications that Fourier reveals. 
The derivative of integer order $n\geq0$ of $u$ with respect to $t$ can be computed as
$$\frac{d^nu}{dt^n}(t)=\frac{1}{2\pi}\int_{-\infty}^\infty dy\,u(y)\int_{-\infty}^\infty
d\omega\frac{d^n}{dt^n}\cos(\omega t-\omega y).$$
Observe that the derivative operator acts only on the cosine in the right hand side.
Similarly, one can represent the integral of $u$ with respect to $t$
with a formula involving the integration in $t$ of the cosine function.
In \cite[Par.~422]{MR1414430}, Fourier
points out that, since $\frac{d^n}{dt^n}\cos(\omega t-\omega y)=\omega^n\cos(\omega t-\omega y+n\frac{\pi}{2})$,
we have
{\small$$\frac{d^nu}{dt^n}(t)=\frac{1}{2\pi}\int_{-\infty}^\infty dy\,u(y)\int_{-\infty}^\infty
d\omega\,\omega^n\cos\Big(\omega t-\omega y+n\frac{\pi}{2}\Big).$$}
Immediately after this equation, he defines fractional derivatives
and integrals: \emph{``The number $n$, that enters in the second member,
will be regarded as any positive or negative quantity. We shall not insist on these applications to
the general analysis,''} see \cite[p.~562]{MR1414430}.
Perhaps because of this last comment, some have regarded Fourier's definition as belonging to the \emph{prehistory} of fractional
calculus, see \cite[p.~xxvii]{MR1347689}. However, Fourier's definitions of fractional derivatives
and integrals are the first ones given for a general function, not just for power functions as Euler and Lacroix did.

The integral representation of $u(t)$ given above is nothing but the
Fourier transform inversion formula. Indeed, by applying the trigonometric identity for
the cosine of a difference of two angles and Euler's formula $e^{it\omega}=\cos(t\omega)+i\sin(t\omega)$,
Fourier's formula reads
\begin{equation}\label{eq:representation}
u(t)=\frac{1}{(2\pi)^{1/2}}\int_{-\infty}^\infty\widehat{u}(\omega)e^{it\omega}\,d\omega
\end{equation}
where the Fourier transform $\widehat{u}$ of $u$ is given by
$$\widehat{u}(\omega)=\frac{1}{(2\pi)^{1/2}}\int_{-\infty}^\infty u(t)e^{-it\omega}\,dt.$$
Obviously, Fourier had a clear understanding of the complex form of his transform
(see, for instance, \cite[Par.~420]{MR1414430}), but he worked mostly with the cosine formulation.
Now, Fourier established that
\begin{equation}\label{eq:derivatives}
\frac{d^nu}{dt^n}(t)=\frac{1}{(2\pi)^{1/2}}\int_{-\infty}^\infty(i\omega)^n\widehat{u}(\omega)e^{it\omega}\,d\omega.
\end{equation}
In other words, differentiation becomes an algebraic operation:
derivatives of integer order $n$ are just multiplication of the Fourier transform $\widehat{u}$
by the homogeneous complex monomial $(i\omega)^n$.
Following Fourier, the number $n$ can now be regarded as any positive or negative quantity.
Thus, Fourier's fractional derivative of $u$ of order $\alpha>0$ is defined by
$$D^\alpha u(t)=\frac{1}{(2\pi)^{1/2}}\int_{-\infty}^\infty(i\omega)^\alpha\widehat{u}(\omega)e^{it\omega}\,d\omega,$$
which is the same as saying that
$$\widehat{D^\alpha u}(\omega)=(i\omega)^\alpha\widehat{u}(\omega).$$
Similarly, the fractional integral of order $\alpha>0$ is
$$\widehat{D^{-\alpha}u}(\omega)=(i\omega)^{-\alpha}\widehat{u}(\omega).$$
It is obvious from the definitions above that the composition formula $D^{\alpha_1}(D^{\alpha_2}u)=D^{\alpha_1+\alpha_2}u$,
the Fundamental Theorem of Fractional Calculus $D^\alpha(D^{-\alpha}u)=D^{-\alpha}(D^{\alpha}u)=u$,
and the consistency limits $\lim_{\alpha\to0}D^\alpha u=u$,
$\lim_{\alpha\to1}D^\alpha u=\frac{du}{dt}$ and, in general, $\lim_{\alpha\to n}D^\alpha u=\frac{d^nu}{dt^n}$ all hold.

Although this seems to work just fine as a good definition of fractional differentiation,
we are now faced with many questions.
For which functions $u$ is $D^\alpha u$ well-defined? How do we perform
the inverse Fourier transform to actually compute $D^\alpha u(t)$ for specific functions $u$?
Does $D^\alpha u(t)$ look like a limit of an incremental quotient, similar to the classical derivative?
For which classes of $u$ does the Fundamental Theorem of Fractional Calculus hold?
In which sense are the consistency limits valid?
But, even before all of that: what is the meaning of the fractional power $(i\omega)^\alpha$
for the purely imaginary complex number $i\omega$? For example, $(2i)^{1/2}=\pm(1+i)$,
so which root should we choose? We can begin to answer these questions by
introducing a powerful tool: the \emph{method of semigroups}.

\section*{The method of semigroups for fractional derivatives}

The definitions of fractional derivative $D^\alpha$
and fractional integral $D^{-\alpha}$ given by Fourier can be realized as the positive
and negative \emph{fractional powers of the derivative operator}, respectively.
Indeed, we do this just by extension of meaning: in the Fourier side,
the derivative is multiplication by $(i\omega)$, the integer power $n$ of the derivative operator is $(i\omega)^n$,
so the fractional power $\alpha$ of the derivative is $(i\omega)^\alpha$.

The method of semigroups is a very general tool to precisely define, characterize,
analyze and use fractional powers of linear operators in concrete problems. In fact, it can be applied to operators such as
the Laplacian  \cite{MR3965397}, the heat operator \cite{MR3709888}, the Laplace--Beltrami operator in a Riemannian
manifold, the discrete Laplacian, the discrete derivative, the wave operator, and many others. It has also been
used in numerical and computational implementation
of fractional operators with finite differences and finite elements methods. The theory was first
established in \cite{MR2754080} for operators on Hilbert spaces, and then extended
in \cite{MR3056307} to operators on Banach spaces. 
We will not describe the whole theory here nor list all of its applications
(to see how it works for the case of the fractional Laplacian
and find more references, we refer to \cite{MR3965397}).
Instead, we will show how the methodology can
unpack Fourier's definition of fractional derivatives and start answering some of the questions left open
at the end of the previous section.
The following discussion is based on \cite{MR3456835}, where the reader can find all the proofs.

From now on, we will just focus on the case $0<\alpha<1$. If we want to analyze higher order fractional derivatives,
say, of order $3/2$ or $\pi$ then, as we can write $D^{3/2}u=D^{1/2}\big(\frac{du}{dt}\big)$
and $D^\pi u=D^{\pi-3}\big(\frac{d^3u}{dt^3}\big)$, we only need to study the operators $D^{1/2}$ and $D^{\pi-3}$.

Let us next recall three important properties of
the Fourier transform. We have already seen two of them: the Fourier inversion formula in \eqref{eq:representation}
and the relation between the Fourier transform and derivatives of order $n$ in \eqref{eq:derivatives}.
The third one is that the Fourier transform of a translation of $u$ corresponds to modulation of $\widehat{u}$
by a complex exponential. Indeed, if the translation operator $T_\tau$ acts on $u$
as $T_\tau u(t)=u(t-\tau)$, then $\widehat{T_\tau u}(\omega)=e^{-\tau(i\omega)}\widehat{u}(\omega)$.
Clearly, we have $T_{\tau_1}(T_{\tau_2}u)=T_{\tau_1+\tau_2}u$
and $\lim_{\tau\to0}T_\tau u=T_0u=u$. These properties imply that the family of operators $\{T_\tau\}_{\tau\geq0}$ is 
a \emph{semigroup}, known as the semigroup of left translations. We call it ``left'' because, for $\tau>0$, $T_\tau u(t)$
looks at the values of $u$ to the left of $t$. The restriction to $\tau>0$ will be apparent soon.

The key to the semigroup method for fractional derivatives
lies in two important integral identities involving the Gamma function $\Gamma$.
By using the Cauchy integral theorem and the unique continuation theorem of complex analysis,
it is proved in \cite[Corollary~2.2]{MR3456835} that, for any $\omega\in\R$ and $0<\alpha<1$,
$$(i\omega)^\alpha=\frac{1}{\Gamma(-\alpha)}\int_0^\infty\big(e^{-\tau(i\omega)}-1\big)\,\frac{d\tau}{\tau^{1+\alpha}}.$$
In particular, this formula implies that $(i\omega)^\alpha$ on the left hand side is chosen from the principal branch
of the multi-valued complex function $z^\alpha$, and this answers
the question of which power should we select.
If we multiply both sides of this identity by $\widehat{u}(\omega)$, the left hand side becomes
$\widehat{D^\alpha u}(\omega)$, while in the integrand on the right hand side we
get $e^{-i\tau\omega}\widehat{u}(\omega)-\widehat{u}(\omega)$, which is the Fourier transform of
$T_\tau u(t)-u(t)=u(t-\tau)-u(t)$ (note that $\tau>0$).
Hence, after inverting back from the Fourier side and making a simple change of variables,
we find the pointwise formula
\begin{equation}\label{eq:Marchaud}
\begin{aligned}
D^\alpha u(t) &= \frac{1}{\Gamma(-\alpha)}\int_0^\infty\frac{T_\tau u(t)-u(t)}{\tau^{1+\alpha}}\,d\tau \\
&= c_{\alpha}\int_{-\infty}^t\frac{u(t)-u(\tau)}{(t-\tau)^{1+\alpha}}\,d\tau.
\end{aligned}
\end{equation}
One impressive aspect of the method of semigroups becomes evident:
we have found the pointwise formula for $D^\alpha u(t)$ without directly computing the inverse
Fourier transform of $(i\omega)^\alpha\widehat{u}(\omega)$
(which should be performed in the sense of tempered distributions because $(i\omega)^\alpha$
is not a bounded Fourier multiplier).

The fractional derivative \eqref{eq:Marchaud} involves a sort of fractional incremental quotient
in which $u(t)$ is compared with $u(\tau)$ through the \emph{interaction} kernel
$1/(t-\tau)^{1+\alpha}$. The kernel becomes singular when $\tau=t$. This gives an idea
that $u$ must have some regularity at $t$ in order to have a well-defined fractional derivative.

Now, to compute $D^\alpha u(t)$ we need to know the values of $u(\tau)$ for \emph{all} $\tau\leq t$.
In other words, $D^\alpha$ is a \emph{nonlocal} operator. Nonlocal also means that
if $u$ has compact support then, in general, $D^\alpha u$ has noncompact support.
This is not the case for the computation of the classical derivative in which it is enough to know $u(\tau)$
for values $\tau$ infinitesimally close to $t$. Furthermore, the support of $\frac{du}{dt}$ is always
contained in the support of $u$. Because of this, we say that classical differential operators are local operators.

Another aspect of \eqref{eq:Marchaud} is that $D^\alpha u(t)$ is \emph{one-sided}
in the sense that we need the values of $u(\tau)$ for  $\tau<t$,
that is, \emph{from the past}.
We should remember that the classical derivative has a hidden \emph{two-sided} structure.
Indeed, in calculus, we first introduce the derivatives from the left and from the right by taking
left- and right-sided limits of incremental quotients, respectively. In the very special
situation in which both limits coincide, we call that common limit \emph{the} derivative.
The difference between left and right derivatives is
very well understood in numerical analysis because of the different properties between the
explicit/forward and the implicit/backwards Euler methods.

As a matter of fact, \cite{MR3456835} shows
that $D^\alpha$ is the fractional power $(D_{\mathrm{left}})^\alpha$ of the left derivative
operator $D_{\mathrm{left}}$. Here $D_{\mathrm{left}}$ is the (negative of the) infinitesimal
generator of the left translation semigroup $T_\tau u(t)=u(t-\tau)$, $\tau>0$,
namely,
$$D_{\mathrm{left}}u(t)=-\lim_{\tau\to0^+}\frac{u(t-\tau)-u(t)}{\tau}.$$
Thus, if $u$ is differentiable, then $D_{\mathrm{left}}u=u'$. One can also
obtain the fractional derivative that looks \emph{into the future} by taking the fractional
power of $D_{\mathrm{right}}$, the derivative from the right,
which is the infinitesimal generator of the right translation semigroup, see \cite{MR3456835}.

The pointwise formula for the fractional integral $D^{-\alpha}u(t)$ can be found in a similar way,
starting with the other important Gamma function integral formula
$$(i\omega)^{-\alpha}=\frac{1}{\Gamma(\alpha)}\int_0^\infty e^{-\tau(i\omega)}\,\frac{d\tau}{\tau^{1-\alpha}}.$$
Parallel as before, it follows that
\begin{equation}\label{eq:Weyl}
\begin{aligned}
D^{-\alpha}u(t) &= \frac{1}{\Gamma(\alpha)}\int_0^\infty\frac{T_\tau u(t)}{\tau^{1-\alpha}}\,d\tau \\
&= c_{-\alpha}\int_{-\infty}^t\frac{u(\tau)}{(t-\tau)^{1-\alpha}}\,d\tau,
\end{aligned}
\end{equation}
which is the negative fractional power of the left derivative operator $(D_{\mathrm{left}})^{-\alpha}$.
For these details, see \cite{MR3456835}.

The left-sided fractional derivative \eqref{eq:Marchaud}, which from now on will be denoted by $(D_{\mathrm{left}})^\alpha$,
is known as the Marchaud--Weyl fractional derivative.
Andr\'e Marchaud found this formula in his 1927 PhD dissertation \cite{MR3532941}.
His derivation, though, followed completely different motivations and arguments.
Indeed, Marchaud wanted to extend the Riemann--Liouville fractional derivative
(which we will not discuss here) to unbounded intervals. Instead, our approach is based on Fourier's original definition \cite{MR1414430}
and semigroups as in \cite{MR3456835}. On the other hand, the fractional integral \eqref{eq:Weyl},
which will be denoted by $(D_{\mathrm{left}})^{-\alpha}$, is known as the
Weyl fractional integral. It was introduced by Hermann Weyl in his 1917 paper \cite{MR3618577}.
Weyl used Fourier series and even found the formula \eqref{eq:Marchaud} for the inverse of $(D_{\mathrm{left}})^{-\alpha}$.
Again, \eqref{eq:Weyl} looks at the values of $u$ in the past
and, like the classical integration operator, is nonlocal.

Since these fractional operators consider the values of $u$ in the past, it seems reasonable to use
\eqref{eq:Marchaud} or \eqref{eq:Weyl} to account for \emph{memory effects}. In probability language,
non-Markovian processes are those in which the evolution of a system
depends not only on the present but also on the past history. As we will see, this intuition amounts for effective models for population growth,
for viscoelastic materials response in mechanics and bioengineering, and for anomalous diffusions
in physics, fluid mechanics and biology.
But before considering these applications, let us outline some elements of the recently developed
analytical theory of (left-sided) fractional derivatives and integrals.
The theory for right-sided fractional calculus can be established without any problems in an analogous way.

\section*{Theory of left-sided fractional derivatives}

It is easy to check that if $u$ is, say, bounded at $-\infty$, and H\"older continuous of order $\alpha+\epsilon<1$
at $t$ from the left, for some $\epsilon>0$, then $(D_{\mathrm{left}})^\alpha u(t)$ is well-defined. More generally, if $u\in C^\beta$
for $\beta>0$ with appropriate decay at $-\infty$, then
$(D_{\mathrm{left}})^\alpha u\in C^{\beta-\alpha}$ and $(D_{\mathrm{left}})^{-\alpha} u\in C^{\beta+\alpha}$. Thus, at the scale
of H\"older spaces, $(D_{\mathrm{left}})^\alpha$ and $(D_{\mathrm{left}})^{-\alpha}$ behave as differentiation
and integration of fractional order $\alpha$, respectively.

Clearly, the fractional derivative of a constant is zero: $(D_{\mathrm{left}})^\alpha1\equiv0$ on $\R$.
By using the definition of the Beta function, it can be seen that, for any $\beta>\alpha$,
\begin{equation}\label{eq:powerfunction}
(D_{\mathrm{left}})^\alpha[(t_+)^\beta]=\frac{\beta\Gamma(\beta)}{\Gamma(1+\beta-\alpha)}(t_+)^{\beta-\alpha}.
\end{equation}
We can also see the influence of the past by modifying the previous function on $(-\infty,0]$.
Indeed, if we define $u_\beta(t)=t^\beta$ for $t>0$ and $u_\beta(t)=1$ (instead of $0$) for $t\leq0$, then
\begin{multline}\label{eq:powermodified}
(D_{\mathrm{left}})^\alpha u_\beta(t)\\
=\frac{\beta\Gamma(\beta)}{\Gamma(1+\beta-\alpha)}(t_+)^{\beta-\alpha}
+\frac{(t_+)^{\beta-\alpha}-(t_+)^{-\alpha}}{\Gamma(1-\alpha)}.
\end{multline}
With a simple change of variables one can also check that, for any $\lambda>0$,
$$(D_{\mathrm{left}})^\alpha(e^{\lambda t})=\lambda^\alpha e^{\lambda t}.$$
From here, we can deduce the fractional derivatives of sine and cosine.
In all these examples, when $\alpha\to1$ we obtain the left derivative of the functions.

Nevertheless, some paradoxes arise in the limit as $\alpha\to0$. Obviously, 
we cannot recover the constant function $1$ from the limit $\lim_{\alpha\to0}(D_{\mathrm{left}})^\alpha1=0$.
More surprisingly, in the limit as $\alpha\to0$ in \eqref{eq:powerfunction} we get $(t_+)^\beta$ back,
but that is not the case of \eqref{eq:powermodified}, in which for $t>0$ we obtain $2t^\beta-1\neq u_\beta(t)$.
This is a consequence of the nonlocality: $u_\beta$ has a \emph{fat tail} at $-\infty$, and this influences
its fractional derivative at $t>0$ more than the \emph{skinny tail} of $(t_+)^\beta$ does.

For sufficiently good functions $u$ and $\varphi$, one can check that
$$\int_{-\infty}^\infty\big[(D_{\mathrm{left}})^\alpha u\big]\varphi\,dt=\int_{-\infty}^\infty u\big[(D_{\mathrm{right}})^\alpha\varphi\big]\,dt.$$
In view of this relation, we can define the fractional derivative of a distribution $u$ as $[(D_{\mathrm{left}})^\alpha u](\varphi)
=u((D_{\mathrm{right}})^\alpha\varphi)$, for suitable test functions $\varphi$. The distributional space
must reflect the one-sided nature of fractional derivatives. It was shown in \cite{MR4062978} that the appropriate
test functions must be supported on intervals of the form $(-\infty,A]$, so as to look at
$(D_{\mathrm{left}})^\alpha u$ \emph{from the left}. Since $(D_{\mathrm{right}})^\alpha\varphi$ will also
have support in $(-\infty,A]$, in the action $u((D_{\mathrm{right}})^\alpha\varphi)$ the only values
of $u$ involved are those \emph{to the left}.

Notice that, if $\varphi$ is smooth, then 
$D_{\mathrm{left}}u(\varphi)=u(D_{\mathrm{right}}\varphi)=-u(\varphi')$.
In other words, due to the regularity of the test functions $\varphi$, in the distributional and weak senses the left derivative
coincides with the classical derivative. Thus, to define left-sided Sobolev spaces, one needs to introduce
weighted spaces that are capable of encoding the underlying left-sided structure. This is accomplished by using one-sided
Sawyer weights $w(t)\in A_p^-(\R)$, which are the good weights for the left-sided Hardy--Littlewood
maximal function
$$M^-u(t)=\sup_{h>0}\frac{1}{h}\int_{t-h}^t|u(\tau)|\,d\tau.$$
(It is important to remark that
this is the original definition of maximal function given by Hardy and Littlewood in \cite{MR1555303}.)
Indeed, $M^-$ is bounded on the weighted space $L^p(\R,w)$, $1<p<\infty$,
if and only if $w\in A_p^-(\R)$, see \cite{MR849466}.
The class $A_p^-(\R)$ is larger than the usual Muckenhoupt class $A_p(\R)$
as the example $w(t)=e^{-t}$ shows. The left-sided weighted Sobolev space
$W^{1,p}(\R,w)$ for the left derivative
is then defined as the set of functions $u\in L^p(\R,w)$ such that $u'\in L^p(\R,w)$, where $w\in A_p^-(\R)$
and $1\leq p<\infty$. The space is consistent with the left-sided fractional calculus
in the sense that if $u\in L^p(\R,w)$ then $(D_{\mathrm{left}})^\alpha u$ is well-defined in the sense
of distributions. Moreover, if $\lim_{\alpha\to1}(D_{\mathrm{left}})^\alpha u=v$ exists in $L^p(\R,w)$
then $u\in W^{1,p}(\R,w)$ and the limit $v$ is exactly $u'$. Conversely, if $u\in W^{1,p}(\R,w)$
then $\lim_{\alpha\to1}(D_{\mathrm{left}})^\alpha u=u'$ in $L^p(\R,w)$ and almost everywhere.
It can also be proved that if $u\in L^p(\R,w)$ and $\lim_{\alpha\to0}(D_{\mathrm{left}})^\alpha u=v$ exists in $L^p(\R,w)$,
then the limit $v$ is indeed $u$ almost everywhere. See \cite{MR4062978}.

The one-sided fractional Sobolev spaces can also be defined.
If $w\in A_p^-(\R)$, $1<p<\infty$, then $W^{\alpha,p}(\R,\omega^p)$ is the set
of functions $u$ for which there is $f\in L^p(\R,w^p)$ such that $u=(D_{\mathrm{left}})^{-\alpha}f$.
In fact, $u$ is in this potential space if and only if $(D_{\mathrm{left}})^\alpha u$
exists in $L^p(\R,w^p)$, see \cite{MR3456835}.

The Fundamental Theorem of Fractional Calculus is established in one-sided weighted $L^p$
spaces. In \cite{MR3456835} it is shown that $u(t)=(D_{\mathrm{left}})^\alpha(D_{\mathrm{left}})^{-\alpha}u(t)$
in the sense of one-sided weighted $L^p$ spaces and almost everywhere. Observe that this result,
which is obvious for $u$ smooth with compact support through the Fourier transform, is not trivial at
all in $L^p$ or in the almost everywhere sense. Indeed, the proof involves intricate integrations
and delicate maximal function estimates.

From the PDE perspective, the nonlocality of $(D_\mathrm{left})^\alpha$ creates an obstacle
in that classical localization techniques of multiplying by a compactly supported test function
and integrating by parts are not directly applicable. This happens because the fractional derivative of a compactly
supported function is not necessarily of compact support. The so-called \emph{extension problem}
gives a local characterization of fractional derivatives in terms of a PDE that involves only classical derivatives.
This extension technique in the PDE context
was first introduced for the fractional Laplacian by Caffarelli and Silvestre (see \cite{MR3965397}
for an overview). The method of semigroups was developed
and used to generalize the extension problem to fractional power operators in \cite{MR2754080,MR3056307}.
For the case of fractional derivatives,
it was shown in \cite{MR3456835} that, if $U(t,y):\R\times[0,\infty)\to\R$ is the solution to
$$\begin{cases}
-D_{\mathrm{left}}U+\frac{1-2\alpha}{y}U_y+U_{yy}=0,&\hbox{for}~t\in\R,~y>0\\
U(t,0)=u(t)&\hbox{for}~t\in\R
\end{cases}$$
then, for some explicit constant $d_\alpha>0$,
$$-\lim_{y\to0^+}y^{1-2\alpha}U_y(t,y)=d_\alpha(D_{\mathrm{left}})^\alpha u(t).$$
The limit can be taken in the classical sense, but also in the sense of distributions, in one-sided weighted
$L^p$ sense and almost everywhere, depending on the regularity of $u$. With the extension,
problems involving fractional derivatives of $u$ can be converted into equivalent local PDE problems for $U$.
The prices to pay, though, are the increase of dimension (passing from a one-dimensional problem for $u$
to a two-dimensional problem for $U$) and the degeneracy of the PDE for $U$.
The advantage is that many analytical and numerical PDE tools can be implemented
in the extension problem in order to obtain estimates and properties of $u$. For instance,
with the extension problem one can prove
Harnack inequalities for nonnegative solutions to $(D_{\mathrm{left}})^\alpha u=0$ in an interval $I\subset\R$,
see \cite{MR3543684} and \cite{MR3709888}.

After this account on the theory of left-sided fractional derivatives, and before presenting
the three applications, we need to address the elephants in the room.

\section*{The elephants in the room}

\begin{figure}
\begin{center}
\includegraphics[width=0.45\textwidth]{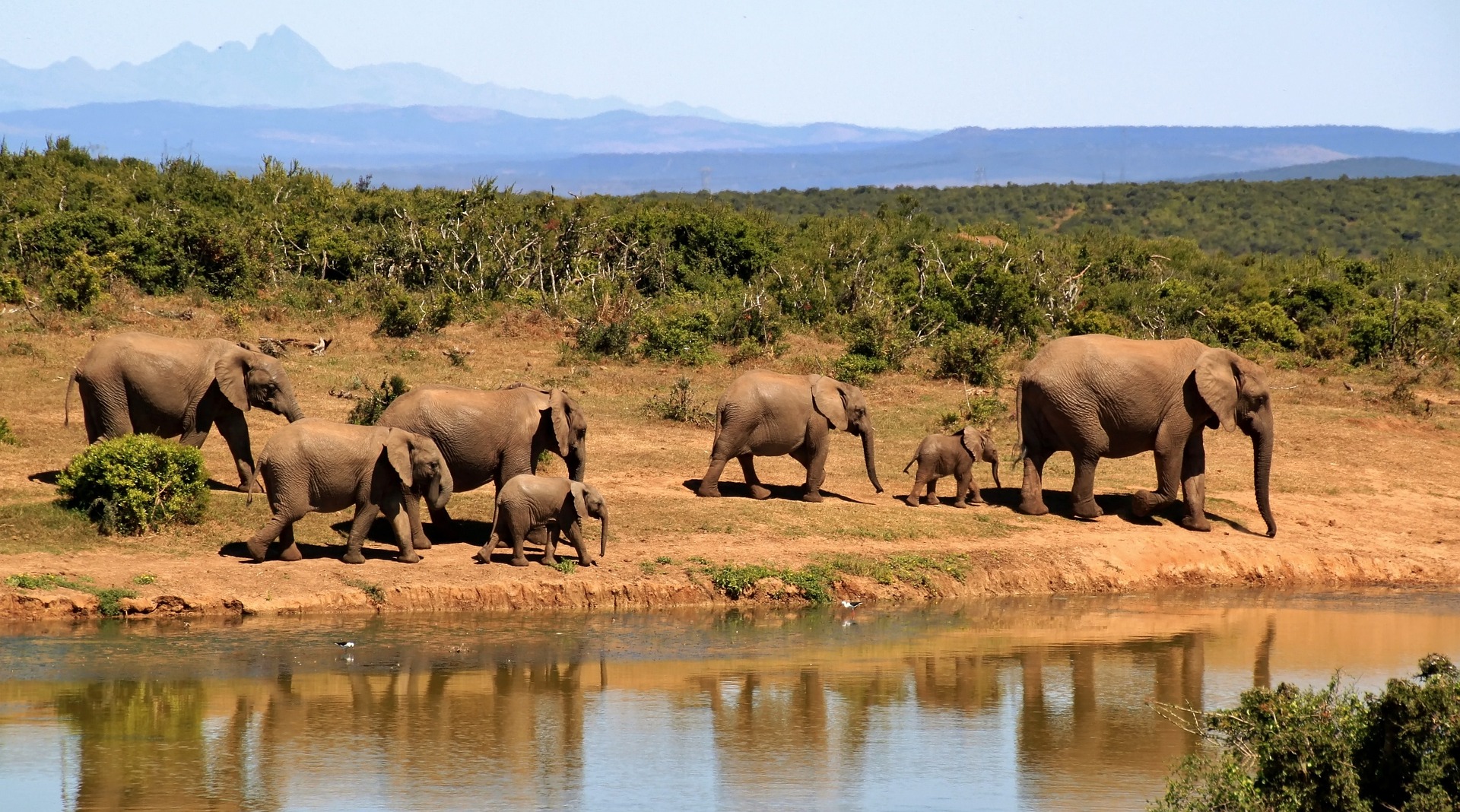}
\caption{The elephants in the room: because of their memory, elephants' travel pattern
follows an anomalous diffusion process. The elephant random walk was introduced in \cite{Schtz2004ElephantsCA}
to model memory effects.}\label{figure:elephants}
\end{center}
\end{figure}

When talking about fractional calculus, it is important to recognize and reflect upon some
of the questions we may encounter.

One could argue that Fourier's fractional derivative is just a mathematical curiosity
that can be studied from a merely academic point of view. A couple of valid questions are:
\emph{what are the applications of fractional derivatives?} and \emph{why do we need fractional derivatives?}
To address these questions, this article will show three applications:
population growth (related to data interpolation), viscoelastic materials (based on Boltzman's
principle of superposition) and anomalous diffusion
(based on observed experiments). These, of course, are not the only places where fractional
calculus arises naturally. The reader is invited to make an effort spending some time learning about other interesting models
in science and engineering by browsing the web, which is filled with survey articles written by
engineers, physicists, biologists and mathematicians,
or by going through the lists collected in \cite{SUN2018213,MR1658022}.

The second elephant in the room, intimately connected to the first one, is the variety of options:
\emph{which fractional derivatives should we use?}
Indeed, there is a zoo of definitions of fractional derivatives out there, such as the ones named after
Riemann, Liouville, Chapman, Caputo, Gr\"unwald, Letnikov and Jumarie, among many others. 
In our opinion, the fractional derivative should be chosen directly from modeling considerations.
We will not address these other definitions nor draw comparisons or connections among them
(a good source is, for instance, \cite{MR1347689}, and many other articles, undergraduate and master's dissertations,
and books that compare them all). Here, we will only focus on Fourier's concept.

Another objection may be that of computational cost:
\emph{it is cheaper to numerically
compute solutions to local PDEs than solutions to nonlocal, fractional order equations.}
We need to face the fact that, at the end of the day, we cannot escape reality:
in many instances, PDEs have proved to be incapable of
reproducing and predicting the observed experiments, but instead nonlocal
fractional models have shown to be the most adequate tool.
Numerical simulation of nonlocal models is a challenging topic in
which much more must be understood. Such active research area is
growing hand in hand with the technological advances that are rapidly decreasing computational costs.

At the same time, it is critical to ask
\emph{why are these elephants in the room?} Today, nobody needs to make a case to justify
why the heat equation is important. One way of explaining this is by reflecting about our training.
Calculus courses have ingrained into our way of thinking \emph{the} concept of instantaneous rate of change,
in which the change is independent of all the previous history (Markovian processes).
We \emph{naturally} think about evolution equation models in terms of classical time derivatives.
In contrast, many systems exhibiting anomalous diffusions and memory effects,
like the way elephants and other wild animals travel, see Figure \ref{figure:elephants}, are not presented in typical undergraduate studies.
Usual undergraduate probability and PDE textbooks
do not develop the idea of modeling memory effects in random walks, where particles
can get stuck in a location for a random period of time, or where plasmas and elephants can travel
long distances without following Fourier's law of diffusion. Instead,
one of the central techniques that is emphasized, and rightly so, is Fourier's method of solving
the heat equation by separation of variables.

Fourier formulated his set of ideas 
in 1807 in his memoir \emph{On the Propagation of Heat in Solid Bodies}.
Highly influential scientists at the time such as Lagrange, Laplace, Monge, Lacroix and Poisson
raised objections to Fourier's trigonometric expansions, see \cite{Her75}. Their opinions prevented the
publication of Fourier's memoir by the French Acad\'emie des Sciences.
It took 15 years for the French scientific community to allow the publication in 1822 of \cite{MR1414430}.
Nowadays, Fourier's work and Fourier series are highly regarded in the scientific community,
but at the time the skepticism of the experts caused controversy, shaking the \emph{status quo}.

It looks like Fourier's ideas on fractional calculus are raising new paradoxes again,
fulfilling the prophetic words of Leibniz.
Who knows, maybe it will take another 15 years for his concept of fractional derivative
to enter into mainstream undergraduate calculus, probability and PDE textbooks and courses.

Let us continue with three applications of Fourier's definition of fractional derivative.

\section*{Population growth with memory}

The most simple ODE model taught in undergraduate calculus
is that of unlimited population growth. Under various simplifying conditions,
the 1798 Malthusian law of of population growth 
establishes that the rate of change $\frac{du}{dt}(t)$ of a population density $u(t)$ at time $t$ is proportional to the current number of
individuals, that is, $\frac{du}{dt}=\lambda u$. Here $\lambda$ is the birth (if positive) or death (if negative) rate.
If the initial population is $u(0)=C>0$, then the unique solution is $u(t)=Ce^{\lambda t}$, for $t\geq0$.
This model predicts that if $\lambda>0$ then there will be such an exponential increase in population that soon there will
be not enough resources on Earth to feed everyone (although, on the other hand, the more people, the more
minds to find solutions to humanity's problems). Obviously, the model
is a good first approximation to cases where the simplifying assumptions are
reasonably met, like radioactive decay or cell reproduction. But it does not directly apply to human population growth.

Populations are influenced by a multitude of
factors such as pandemics, migration, cultural changes, natural disasters, social media, etc.
In other words, populations have \emph{memory}.
Since left-sided fractional derivatives take into account memory effects,
one can propose a similar model of population growth with memory:
$(D_{\mathrm{left}})^\alpha u(t)=\lambda u(t)$, for $t>0$.
The correct \emph{initial} condition is the form $u(t)=C(t)$, for all $t\leq0$.
That is, we need to know the historic population $C(t)$ until the initial time $t=0$. In \cite{MR3557159},
the authors use the World Population data from 
the United Nations to adjust the value of $\alpha$ that best fits the data, obtaining a much better error of
approximation for the fractional model than the classical ODE model.
This is a very simple instance where the classical local equation cannot fit the data, but the nonlocal one does
so effectively. In addition, \cite{MR3557159} analyzes other blood alcohol level and video tape problems.

Notice that the classical exponential solution $e^{\lambda t}$ is in fact an eigenfunction for the derivative operator.
The fractional counterpart of the exponential function is 
the Mittag-Leffler function, which is defined by the power series
$$E_{\alpha,\beta}(t)=\sum_{k=0}^\infty\frac{t^k}{\Gamma(\alpha k+\beta)}.$$
For $\alpha,\beta>0$, it has an infinite radius of convergence.
In particular, $E_{1,1}(\lambda t)=e^{\lambda t}$. If we define
$$u(t)=E_{\alpha,1}(\lambda(t_+)^\alpha)$$
then, by using \eqref{eq:powerfunction}, the reader can verify that
$$(D_{\mathrm{left}})^\alpha u=\lambda u.$$
Thus, the Mittag-Leffler functions are the eigenfunctions of the fractional derivative.
Because of this, we say that they are the \emph{fractional exponentials}.

\section*{Viscoelasticity: materials with memory and Boltzmann}

In classical continuum mechanics, elastic materials are those that return to their original shape
after applied forces have been removed, like steel or concrete under small deformations.
Hooke's law for an ideal elastic solid establishes that the stress $\sigma$ (the internal force in a material per unit area)
is proportional to the strain $\varepsilon$ (the deformation or elongation with respect to the original length), that is,
$\sigma=E\varepsilon$. In this one-dimensional model, $E$ is a material constant, called the elastic or Young modulus,
that can be experimentally measured.
Hooke's equation is pictorially represented by a spring.
On the other hand, liquids, gasses and plasmas deform when subjected to a force.
Newton's law for an ideal fluid says that stress is proportional to the velocity of the deformation,
namely, $\sigma=\eta\frac{d\varepsilon}{dt}$. Here $\eta$ is the fluid dependent,
experimentally observed viscosity coefficient.
This equation is visually represented by a dashpot.

In real life, there are many materials that exhibit both elastic and fluid characteristics
such as wood, asphalt, baker's dough, lead wires, certain polymers, rubber,
clay, gels, metals near melting temperature and even biological tissues like skin.
If we stretch a rubber band for some time, it will lose strength and not return
to its original shape. Similarly, wooden shelves in a library slowly deform due to the weight
of books and will not return to their original shape after the load has been removed but remain warped.
As we age, our skin looses its strength and we develop wrinkles that we might begin to fight using various lotions
or surgical procedures.

In some sense, viscoelastic
materials are in between elastic materials (they do not return to their original shape) and fluids
(they do not continue deforming indefinitely).

Even more in contrast, the examples above show that, unlike elastic solids or fluids, viscoelastic
materials have \emph{memory}: the deformation and the internal forces react by accumulating
the history of all the applied loads and strains.
More precisely, what is particularly interesting about viscoelastic materials is that they have \emph{fading} memory,
see \cite{MR0158605}. Intuitively speaking, ideal liquids ``instantly forget'' where they were, and the Navier--Stokes
equations describing them involve only local time derivatives of the velocity field. Ideal elastic materials have a ``perfect memory''
in that they always remember where they started by returning to their original
configuration, and so have ``no memory'' of previous forces that have since been released.
Viscoelastic materials are intermediate: they remember where they were recently but,
as the internal molecular structure changes irreversibly, forget where they started.
The so-called \emph{decay theory} in psychology and neuroscience tries to explain
the process and causes of fading memory observed in humans.

Early viscoelastic models were created by using
combinations of springs (elastic) and dashpots (viscous) connected in series or in parallel.
Nevertheless, these models have shown to be insufficient to describe the complex behavior of most
viscoelastic materials. Furthermore, it was observed that one would need a very large amount
of springs and dashpots to accurately approximate some materials responses, leading to large systems of
higher order differential equations.

Since fractional derivatives interpolate between the identity operator and the usual derivative,
an intermediate model between springs and dashpots can be given as
$\sigma=k_\alpha D^\alpha\varepsilon$, for some fractional derivative $D^\alpha$ of order $0<\alpha<1$
and some material constant $k_\alpha$.
Mathematically speaking, this is a reasonable model 
between $\alpha=0$ (elastic, $k_0=E$) and $\alpha=1$ (viscous, $k_1=\eta$).
But, after recalling the elephants in the room, we still need a justification of this model coming from physical principles.

\begin{figure}
\begin{center}
\includegraphics[width=0.45\textwidth]{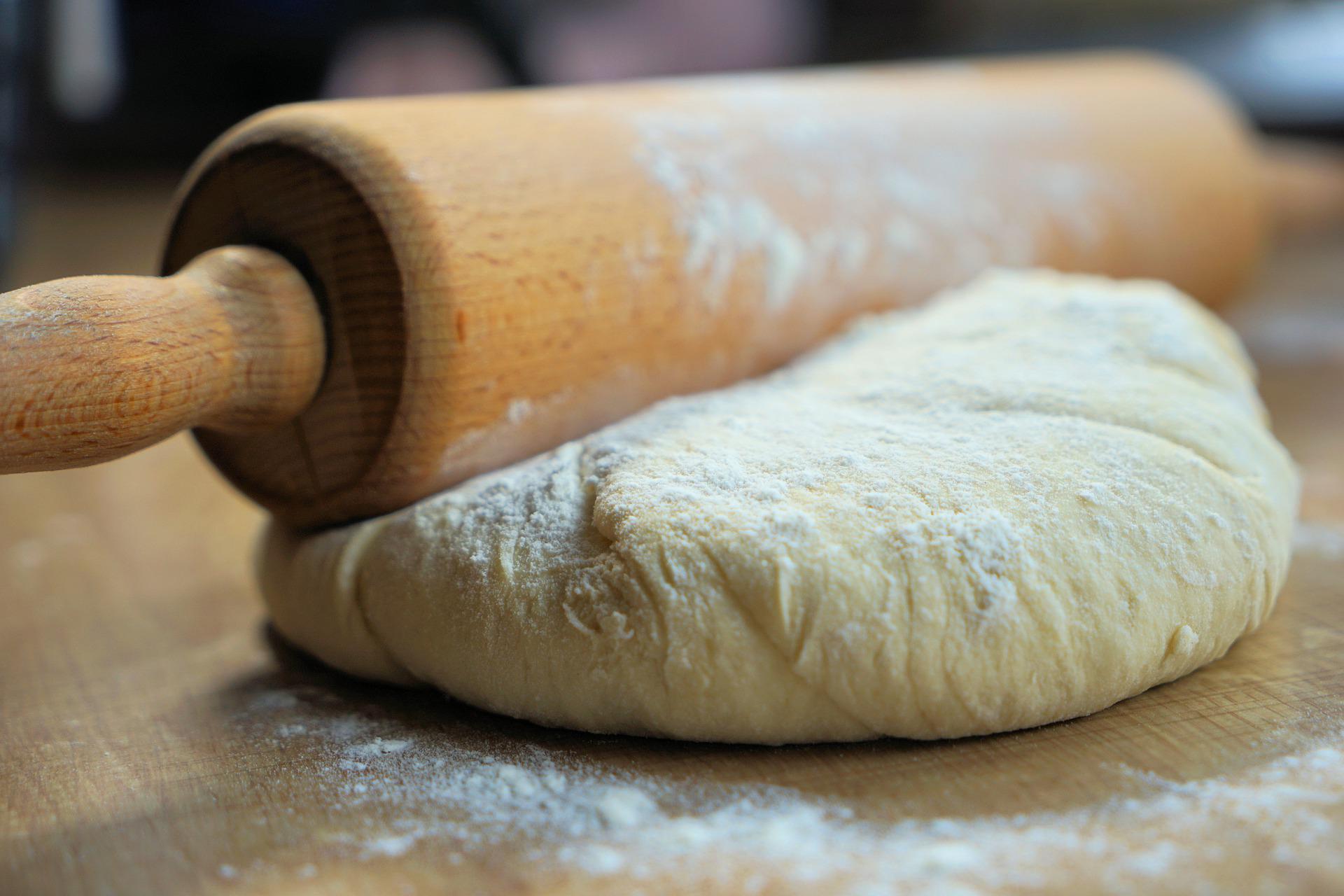}
\caption{Pizza dough is a viscoelastic material with stress relaxation modulus $G(t)=ct^{-0.36}$.}\label{figure:dough}
\end{center}
\end{figure}

To experimentally study the behavior of a viscoelastic material,
one can conduct various tests. The stress relaxation test
measures stress within the material following a fixed displacement that is kept constant with time.
Thus, we apply a step in strain $\varepsilon(t)=\varepsilon_0H(t)$, where $H(t)$ is the Heaviside function and $t$ is time, and calculate
the time-dependent stress $\sigma(t)$. In linear materials, the stress is proportional to the strain, so
$\sigma(t)=\varepsilon_0G(t)$, where $G(t)$ is the so-called stress relaxation modulus.
The relaxation modulus $G(t)$, which is a measure of the material's fading memory,
can be found experimentally. For example, for flour dough, $G(t)=kt^{-0.36}$
\cite[p.~278]{magin2006fractional} (we are neglecting units), see Figure \ref{figure:dough}.
In fact, for many viscoelastic materials, $G(t)$ is a constant multiple of $t^{-\alpha}$, for some $0<\alpha<1$.

In 1874, Ludwig Boltzmann proposed in \cite{Boltzmann}
his principle of superposition to account for memory effects in strained materials.
The principle states that, provided that there is a linear relation between stress and strain, the stress produced by any number
of applied strains is the sum of the stress produced by each of the individual strains when acting alone.
Recall that, in the relaxation test, $\sigma(t)=\varepsilon(0)G(t)$, where $G(t)$
is the stress at time $t$ owing to a unit strain increment of size $\varepsilon(0)$ at time $t=0$.
By the principle of superposition, with another strain increment at time $\Delta\tau>0$,
the stress becomes
$$\sigma(t)=\varepsilon(0)G(t)+\big(\varepsilon(\Delta\tau)-\varepsilon(0)\big)G(t-\Delta\tau).$$
Then, after $N$ increments,
\begin{align*}
\sigma(t)&=\varepsilon(0)G(t)\\
&\quad+\sum_{n=1}^N\frac{\varepsilon(n\Delta\tau)-\varepsilon((n-1)\Delta\tau)}
{n\Delta\tau-(n-1)\Delta\tau}G(t-n\Delta\tau)\Delta\tau.
\end{align*}
In the limit as $\Delta\tau\to0$,
$$\sigma(t)=\varepsilon(0)G(t)+\int_0^t\varepsilon'(\tau)G(t-\tau)\,d\tau.$$
As we said before, we know from experiments that $G(t)=kt^{-\alpha}$, for some $0<\alpha<1$.
A more flexible model that does not involve the derivative of strain in the equation and accounts for all the history of the material
(not necessarily starting at time $t=0$) can be obtained as follows.
We extend $\varepsilon(\tau)\equiv\varepsilon(0)$ for all past times $\tau<0$, extend the integral
from the interval $(0,t)$ to $(-\infty,t)$ (notice $\varepsilon'(\tau)=0$ for $\tau<0$), write
$\varepsilon'(\tau)=\frac{d}{d\tau}\big(\varepsilon(\tau)-\varepsilon(t)\big)$ in the integrand and integrate
by parts to obtain
\begin{align*}
\sigma(t)&=\varepsilon(0)G(t)+\int_{-\infty}^t\big(\varepsilon(\tau)-\varepsilon(t))G'(t-\tau)\,d\tau \\
&=\varepsilon(0)G(t)+c_{k,\alpha}(D_{\mathrm{left}})^\alpha\varepsilon(t),
\end{align*}
where $c_{k,\alpha}$ depends on $k$ and $\alpha$. The decay of the kernel in the fractional
derivative model accounts for the material's fading memory.

A feature of the fractional derivative model is that it provides an adjustable ``material memory''
parameter $\alpha$ for describing the stress/strain behavior of viscoelastic materials that can be experimentally measured.

We mentioned before that features of viscoelastic materials appeared to be better captured by
special combinations of large numbers of springs and dashpots.
One can ask the question: \emph{what combinations of springs and dashpots can give rise to a fractional
derivative model?} It has been shown that hierarchical arrangements of springs and dashpots,
such as infinite ladders, trees and fractal networks,
do in fact obey the fractional viscoelastic constitutive equation in the limit.

Viscoelasticity is a fascinating area of continuum mechanics and bioengineering that the reader
is invited to explore, for instance, in \cite{magin2006fractional}, see also \cite[Chapter~10]{MR1658022}
and references therein.

\section*{Anomalous diffusion is normal}

In the heat equation for an ideal one-dimensional metal rod
 $u_t=\frac{k}{2}u_{xx}$, the function $u(x,t)$ denotes the absolute temperature at the point $x$ at time $t$,
and $k/2>0$ is the diffusivity constant. In deriving this equation, Fourier's main observation
was that the flow of heat between two adjacent molecules is proportional to the
extremely small difference of their temperatures. In other words, the flow of caloric energy
is from regions of higher concentration to regions of lower concentration,
an intuitive, experimentally observed fact. Paradoxically, this is a parallel process
to that of the random movement of a pollen particle in suspension on the surface of water
as observed under the microscope by Scottish botanist Robert Brown in 1927.
Einstein, in his study of Brownian motion, derived the heat equation from first principles
assuming that the direction of motion of the particle is ``forgotten''
after an infinitesimally small period of time.

Let us derive the heat equation by considering a one-dimensional random walk.
We choose a small step size $\Delta x>0$ and a small time interval $\Delta\tau>0$.
Consider a random walker that moves randomly along the $x$ axis
according to the following rules. During the interval of time $\Delta\tau$, the walker takes one step
of size $\Delta x$, starting from, say, $x=0$. The walker moves either to the left or to the right with probability $1/2$,
independently of the previous steps. We would like to compute the probability $u(x,t)$ of finding the walker at position $x$
at time $t$. The process has no memory since each step is independent from the previous ones.
The law of total probability gives that 
$$u(x,t)=\frac{1}{2}u(x-\Delta x,t-\Delta\tau)+\frac{1}{2}u(x+\Delta x,t-\Delta\tau).$$
Indeed, at time $t$, the walker may arrive to $x$ from either previous positions $x-\Delta x$ or
$x+\Delta x$ with probabilities $1/2$.
If we denote the second order incremental quotient of $u$ in the space variable by
$$\delta_2u(x,t)=\frac{u(x-\Delta x,t)+u(x+\Delta x,t)-2u(x,t)}{(\Delta x)^2},$$
then, after subtracting $u(x,t-\Delta\tau)$ to both sides, the identity above becomes
$$\frac{u(x,t)-u(x,t-\Delta\tau)}{\Delta\tau}=\frac{(\Delta x)^2}{2\Delta\tau}\delta_2u(x,t-\Delta\tau).$$
In the limit $\Delta x,\Delta\tau\to 0$, assuming that $\frac{(\Delta x)^2}{\Delta\tau}\to k$, we arrive to the heat equation
$$D_{\mathrm{left}}u=\tfrac{k}{2}u_{xx}.$$
To find the fundamental solution, suppose that we are given the initial condition $u(x,0)=\delta_0(x)$,
a Dirac delta or \emph{unit impulse} concentrated at the origin.
Then, applying the Fourier transform in space $\widehat{u}(\omega,t)$, we find a family of ODEs
parametrized by the Fourier variable $\omega\in\R$:
$$\begin{cases}
D_{\mathrm{left}}\widehat{u}=-\tfrac{k}{2}|\omega|^2\widehat{u},&t>0\\
\widehat{u}=1&t=0.
\end{cases}$$
The solution is $\widehat{u}(\omega,t)=e^{-t\frac{k}{2}|\omega|^2}$. After inverting
the Fourier transform, we find that
$$u(x,t)=\frac{1}{(2\pi kt)^{1/2}}e^{-|x|^2/(2kt)}.$$
This is the classical Gaussian or normal distribution for normal diffusion,
giving the probability of finding the (memoryless) random walker at position $x$ at time $t$.
The mean or average is $\mu=0$
and the standard deviation is $\sigma=(kt)^{1/2}$. In particular, as time goes by, the walker deviates from the origin
an average distance $(kt)^{1/2}$.
The mean square displacement or second moment
$\langle x^2\rangle$ is $kt$. In fact, in the limit we took above, the scaling $(\Delta x)^2=k\Delta\tau$
says that the mean square displacement is proportional to the waiting time $\Delta\tau$ between steps.

There are many instances where the distribution of a quantity is not \emph{normal}.
For instance, the distribution of wealth is not distributed according to a Gaussian. Indeed,
the spread between extreme poverty and wealth is so large
that talking about the average wealth makes no sense.
Such \emph{abnormal} phenomenon follows a Pareto distribution,
in which the probability of finding extremely wealthy people is positive and 
the mean is infinite.

In the mid 1970s, researchers began to pay much more attention to these non-Gaussian processes
and other instances where Einstein's assumptions do not hold.
In \cite{Scher:1975aa}, Scher and Montroll observed in photocopiers and laser printer machines
that the transport of electrons did not follow the diffusion equation. The hypothesis
is that electrons get stuck in ``holes'' within the surface of
amorphous semiconductors for a time and then are released due to a temperature potential.
Physicists refer to this as \emph{diffusion on disordered media}, or simply as \emph{anomalous diffusion}.
It was also observed that the probability distribution of waiting time
in between steps $\psi(\tau)$ is proportional to a Pareto power law $\tau^{-(1+\alpha)}$, for $0<\alpha<1$,
for large waiting times $\tau$.
Intuitively, if the waiting times between steps is large compared to the step size,
then the random walker will not deviate from its initial position as much as a Gaussian random walker
would do, a process that is known as \emph{subdiffusion}.
In these processes, the relation $\langle x^2\rangle=kt$ is lost, but the new subdiffusion power law
$\langle x^2\rangle=k_\alpha t^\alpha$, for some $0<\alpha<1$ and some $k_\alpha>0$, is observed.

Many experiments and observations in nature have shown anomalous diffusions, including the diffusion of lipids and receptors in cell membranes,
the transport of molecules within the cytosol and the nucleus, the travel strategies of wild animals, the sleep-wake
transitions during sleep, the propagation of electric currents on cardiac tissue,
the avalanche-like behavior of plasma particles,
and the fluctuations of the stock market.
In fact, the claim in \cite{Klafter2005AnomalousDS} is that \emph{``the clear
picture that has emerged over the last few decades is that although these phenomena are called anomalous,
they are abundant in everyday life: anomalous is the new normal!''}

New mathematical models to describe
anomalous diffusion have been developed in recent years, including continuous time random walks,
elephant random walks that take into account memory \cite{Schtz2004ElephantsCA}, and nonlocal master equations,
among others.

We do not have the space to enter into more details about the fascinating
world of anomalous diffusions. We invite the reader to explore the literature,
suggesting in particular the popular article \cite{Klafter2005AnomalousDS} that contains many experimental
examples, the detailed surveys \cite{MR1809268} and \cite{MR1937584}, as well as the original
article by Scher and Montroll \cite{Scher:1975aa}.
We will restrict to describing a random walk with memory effects that can easily be
introduced in any undergraduate calculus, probability or differential equations class.

Let us consider then a random walker that follows the same space dynamics as before,
moving a step of size $\Delta x>0$ either to the left or to the right with probability $1/2$,
but also stops at each location for a random period of time. Hence,
there is a waiting time in between steps that is random as well.
In other words, the walker undergoes \emph{memory},
since the next step to the left or to the right happens after a random time
drawn from a distribution of waiting times $\psi(\tau)$.
We are interested in the probability $u(x,t)$ of the walker having just arrived at position $x$ at time $t$.
The law of total probability gives
\begin{align*}
u(x,t)&=\sum_{n=1}^\infty\bigg[\frac{1}{2}u(x-\Delta x,t-n\Delta\tau)\\
&\quad\qquad+\frac{1}{2}u(x+\Delta x,t-n\Delta\tau)\bigg]\psi(n).
\end{align*}
The term in brackets above is related to the probability of arriving at $x$ from either $x-\Delta x$ or $x+\Delta x$,
and those events occur with probability $1/2$. The infinite sum factors in the fact that the walker
could have been at those positions not only at the previous time $t-\Delta\tau$, but may had been sitting there for a
period of time $t-n\Delta\tau$, with the probability of a waiting time of length $n\Delta\tau$
being $\psi(n)$. In accordance with many of the experimental observations
mentioned before, we now assume that $\psi(n)=d_\alpha n^{-(1+\alpha)}$, for some
$0<\alpha<1$, where $d_\alpha>0$ is chosen so that $\sum_{n=1}^\infty\psi(n)=1$.
Using this, we can write the equation above in terms of second order incremental quotients in space as
\begin{multline*}
\sum_{n=1}^\infty\frac{u(x,t)-u(x,t-n\Delta\tau)}{(n\Delta\tau)^{1+\alpha}}(\Delta\tau)\\
=\frac{(\Delta x)^2}{2d_\alpha(\Delta\tau)^{\alpha}}\sum_{n=1}^\infty\delta_2u(x,t-n\Delta\tau)\psi(n).
\end{multline*}
In the limit $\Delta x,\Delta\tau\to 0$, assuming that $\frac{(\Delta x)^2}{d_\alpha(\Delta\tau)^\alpha}\to k_\alpha|\Gamma(-\alpha)|$,
we arrive to the time-fractional equation
$$(D_{\mathrm{left}})^\alpha u=\tfrac{k_\alpha}{2}u_{xx}.$$
To find the fundamental solution, suppose that we are given the past condition $u(x,t)=\delta_0(x)$,
a unit impulse concentrated at the origin, for all times $t\leq0$.
Then, by applying the Fourier transform in space,
$$\begin{cases}
(D_{\mathrm{left}})^\alpha\widehat{u}=-\tfrac{k_\alpha}{2}|\omega|^2\widehat{u},&t>0\\
\widehat{u}=1&t\leq0.
\end{cases}$$
We have already encountered this problem in the population growth model example.
The solution is given in terms of the Mittag--Leffler function as
$$\widehat{u}(\omega,t)=E_{\alpha,1}(-\tfrac{k_\alpha}{2}|\omega|^2(t_+)^\alpha).$$
Using the scaling properties of the Fourier transform, we find that, for $t>0$,
$$u(x,t)=\frac{1}{(\tfrac{k_\alpha}{2}t^\alpha)^{1/2}}H_\alpha\bigg(\frac{|x|}{(\tfrac{k_\alpha}{2}t^\alpha)^{1/2}}\bigg).$$
The profile function $H_\alpha(r)$, $r>0$, is a Fox--Wright function, see, for example, \cite{MR1809268}
or \cite{MR1937584} for this
particular special function.

\section*{Dedication}

This article is dedicated to the loving memory of my friend Roberto A.~Mac\'ias, an extraordinary mathematician
and exemplary person.

\section*{Acknowledgements}

The author is grateful to Jos\'e L.~Torrea and Mary Vaughan for useful discussions and suggestions that helped
improve the presentation of this article. Mary Vaughan provided \cite{MR0141575} and \cite{magin2006fractional}.

\bibliography{StingaRefs}

\end{document}